\documentclass[12pt,reqno]{amsart}
\usepackage{latexsym, amsmath, amscd, amssymb, amsthm}
\usepackage{bm}
\usepackage{mathrsfs}
\usepackage{xspace}
\usepackage{natbib}
\usepackage{setspace}
\usepackage{enumerate}

\usepackage{verbatim}
\usepackage[english]{babel}
\usepackage{epsfig}
\usepackage{algorithm}
\usepackage[noend]{algpseudocode}
\RequirePackage{eucal}
\RequirePackage{eufrak}
\newtheorem{lemma}{Lemma}[section]
\newtheorem{theorem}[lemma]{Theorem}

\theoremstyle{definition}
\newtheorem{definition}[lemma]{Definition}

\newtheorem{example}[lemma]{Example}
\newtheorem{remark}[lemma]{Remark}

\theoremstyle{remark}
\newcommand{\op}{\operatorname}
\newcommand{\lspace}{\vskip7pt}
\newcommand{\PolRing}{K[x_{1},\ldots,x_{n}]}
\begin{document}
\begin{abstract}
 A motivation to study Gr\"{o}bner theory for fields with valuations
 comes from tropical geometry, for example, they can be used to
 compute tropicalization of varieties \citep{maclagan2009introduction}.
 The computational aspect of this theory was first studied in (Chen \& Maclagan, 2013). In
 this paper, we generalize this Gr\"obner basis theory to free modules over
 polynomial rings over fields with valuation. As the valuation of
 coefficients is also taken into account while defining the initial
 term, we do not necessarily get a monomial order. To overcome this
 problem we have to resort to other techniques like the use of ecart
 function where the codomain is the well-ordered set $\mathbb{N}$, and
 thereby give a method to calculate 
 the Gr\"{o}bner basis for submodules generated by homogeneous
 elements.
 Using this, we show how to compute Hilbert polynomials for graded
 modules.
\end{abstract}

\title{Gr\"{o}bner Basis Theory for Modules over Polynomial Rings over Fields with Valuation}
\author[Aritra Sen and  Ambedkar Dukkipati]{Aritra Sen and  Ambedkar Dukkipati}

\email{a.sen@csa.iisc.ernet.in \\ ad@csa.iisc.ernet.in}
\address{Dept. of Computer Science \& Automation\\Indian Institute of Science, Bangalore - 560012}
\maketitle
\section{Introduction}
\nocite{*}%

The Gr\"{o}bner basis theory of modules has several applications in
constructive module theory, for example, in the computation of syzygies,
annihilator of a module and Hilbert polynomial of a graded module
\citep{kreuzer2005computational,spear1977constructive,schreyer1986syzygies}. It
has also various applications in homological algebra, for example in
the computation of Ext and Tor. In computational algebraic geometry,
it is used to compute minimal free resolutions of graded finite modules.  

On the other hand, recently Gr\"{o}bner basis theory for polynomial
rings that takes valuation of the underlying field into consideration
has been studied by~\citet{chan2013groebner}. 
One motivation for this is its various applications
in tropical algebraic geometry.
The other motivation comes from computational aspects of Gr\"{o}bner
theory. They can lead to Gr\"{o}bner bases that are much smaller than
the standard Gr\"obner basis (Remark 4.5). In~\citep{chan2013groebner}
a normal form algorithm has been presented that leads to an algorithm
to compute Gr\"{o}bner basis in this case. 
In this paper we generalize this theory to free modules over
polynomial rings over fields with valuation.  

\subsection*{Contributions}
In Gr\"{o}bner basis theory for modules over polynomial rings over
fields with valuation, since the definition of order on monomials
involves valuations of coefficients it is not possible to generalize
this to the case of modules. So, to overcome this problem we have to
resort to other techniques like the use of an ecart function  where
the codomain is the well-ordered set $\mathbb{N}$. Using  this, we
derive a Buchberger-like criterion for Gr\"{o}bner basis and hence
an algorithm for computing the Gr\"{o}bner basis. One
advantage of this approach is 
that it can lead to smaller Gr\"obner basis. With standard Gr\"obner
basis the initial submodule generally grows with degree $\delta$. One
particular example where the size of the standard Gr\'{o}bner basis
grows linearly with $\delta$ is presented in
\citep{chan2013groebner}. Here, we give an example of a family of
submodules where the size of initial submodules remain constant. Also,
with a slight modification, we show how these ideas can be  ported to
free modules over the polynomial ring
$\mathbb{Z}/p^{\ell}\mathbb{Z}[x_1,\ldots,x_n]$.   

\subsection*{Organization}
The rest of the paper is organized as follows. In Section~2, 
we present preliminaries on fields with valuation and Gr\"{o}bner
basis over fields with valuation. In Section~3, we introduce free
modules over polynomials rings on fields with valuation and a normal
form algorithm for them. We present a  Buchberger-like criterion for
Gr\"{o}bner basis of submodules using the normal form 
algorithm of the previous section and then present an algorithm to 
compute them in Section~4. In Section~5, we show how one can use the
algorithm of Section~4 to compute the Hilbert function of a graded
module. In Section~6, we introduce  free modules over
$\mathbb{Z}/p^{\ell}\mathbb{Z}[x_1,\ldots,x_n]$ and present a normal form
algorithm for them similar to that of Section~3. A  Buchberger-like
criterion for Gr\"obner basis of submodules  of modules over
$\mathbb{Z}/p^{\ell}\mathbb{Z}[x_1,\ldots,x_n]$ and algorithm to compute
the Gr\"{o}bner basis is presented in Section 6. Finally, we give
concluding remarks in Section~7.  

\section{Background}
Throughout this paper, $K$ denotes a field and $\mathbb{N}$ the set of 
natural numbers including zero. For any positive integer $n$, $[n]$
denotes the set $\{1,\ldots ,n\}$. A polynomial ring in indeterminates
$x_1,\ldots,x_n$ over $K$ is denoted by $K[x_1,\ldots , x_n]$. For any
$\alpha \in \mathbb{N}$ a monomial in indeterminates $x_{1},\ldots,x_{n}$
is written as $x^{\alpha}$ and $|\alpha|$ denotes the sum
$\sum\limits_{i=1}^{n} \alpha_{i}$. An arbitrary polynomial $f \in
K[x_{1},\ldots,x_{n}]$ is written as 
\begin{displaymath}
   f = \sum_{\alpha \in \Lambda} a_{\alpha} x^{\alpha} \enspace, 
\end{displaymath}
where $a_{\alpha} \in K$, $\alpha \in {\mathbb{Z}}_{\geq 0}^{n}$ and 
$\Lambda \subseteq {\mathbb{Z}}_{\geq 0}^{n}$ a finite set,
that is support of the polynomial $f$, denoted by $\op{supp}(f)$.
By monomial we mean $x^{\alpha}$, by term we mean $a_\alpha x^\alpha$.

Let $S$ be a finite a set then $|S|$ denotes the number of elements in
the set. Let $w=(w_1,\ldots,w_n) \in \Bbb{R}^n$ and
$u=(u_1,\ldots,u_n) \in \Bbb{R}^n$ then $w.u$ represents the sum
$\sum\limits_{i=1}^{n} w_iu_i$. Monomial order on
$K[x_{1},\ldots,x_{n}]$ is denoted by $\prec$ and for every polynomial
$f \in K[x_1, \ldots , x_n]$, $\op{in}_{\prec}(f)$ denotes the initial
term with respect to $\prec$. Let $I$ be an ideal in
$K[x_1,\ldots,x_n]$, then $\op{in}_\prec(I)=
\langle \op{in}_\prec({f}): {f} \in I \rangle$. $\mathbb Z_{p^{\ell}}$ represents the finite ring $\mathbb{Z}/p^{\ell}\mathbb{Z}$.

\begin{definition}
A field with valuation is an ordered pair of a field $K$ and function
$v$, $(K,v)$ such that  
\begin{enumerate}
\item $v$ is a group homomorphism from $K^*$ to $(R,+,0)$,
\item $v(a+b) \geq \operatorname{min}\{v(a), v(b)\}$
  for all $a,b \in K^*$, and
\item $v(a)= \infty$ iff $a=0$.
\end{enumerate}
\end{definition}
The image of the valuation map is denoted by $\Gamma$. Let $R_K$ be
the set of all field elements with valuation greater or equal to $0$,
i.e., $R_K= \{ c \in K: v(c)\geq 0 \}$. Then $R$ is a local ring with
maximal ideal $J_K= \{c \in K: v(c) > 0\}$ and  $\Bbbk=R_K/J$ is the
residue field of $K$. Let $a \in R$, then $\overline{a}$ represents
its image in the residue field $\Bbbk$. 
\begin{example}
The most common example of a field with valuation is $\mathbb{Q}$ with
$p$-adic valuation $\op{val}_p(.)$, where $p$ is a prime
number. Let $q \in \mathbb Q$, then $\op{val}_p(q)=c$, where
$q=p^ca/b$, such that $p$ does not divide $a$ and $b$. For example,
$\operatorname{val}_3(15/4)=1,\operatorname{val}_2(5/12)=-2$. 
\end{example}
\begin{example}
Consider the field of Puiseux series  $K\{\{t\}\}$,
which is the algebraic closure of the field of
Laurent series when $\op{char}(K)=0$. 
The map, $\operatorname{val}: K\{\{t\}\} \rightarrow \mathbb R$
that takes a Puiseux series and returns the lowest exponent is
a valuation. For example,let $ f(t)= 3t^{-3} + 6t^{-1}+ \ldots $
then $\op{val}(f(t))= -3$. 
\end{example}

\begin{definition}
Let $f= \sum\limits_{u \in\Lambda} c_{u}x^{u} \in
K[x_1,\ldots,x_n]$, be a polynomial, where $c_{u} \in K$, and
$\Lambda$ is the support of $f$. Let $w \in 
\Gamma^n$, $v$ be a valuation of $K$ and  $W=
\min\limits_{u\in\Lambda} \{v(c_{u})+w.u\}$.
The initial form with respect to to $w$ is defined as 
\begin{displaymath}
\op{in}_w(f) = \sum_{v(c_u)+w.u=W} \overline{ c_ut^{-v(c_u)}}x^u.
\end{displaymath}
\end{definition}
Note that  in the above definition of $W$, `$\operatorname{min}$'
can be replaced by `$\op{max}$' since for every valuation $v$ and
weight vector $w$ one can find a weight vector $w'$ and valuation $v'$
such that both cases coincide.   
\begin{example}
Consider the polynomial $f$ over the field of Puiseux
series. $f=(1+t^2)x+2t^2y+t^3z$. Let $w=(1,1,1)$. The
intial form of $f$ with respect to $w$ is $\overline{(1+t^2)}x=x$. 
\end{example}
\begin{example}
Consider the polynomial $f$ over $\mathbb Q$ with $2$-adic valuation,
$f=2x+5y^2+3xyz$. Let $w=(1,1,1)$. Then initial form of
$\op{in}_w(f)=\overline{1}x+\overline{5}y^2=x+y^2$. 
\end{example}

Given an ideal $I$ in $K[x_{1},\ldots,x_{n}]$ and a weight vector $w
\in \Gamma^{n}$, $\op{in}_w(I)$ denotes the ideal $\langle \op{in}_w({f}): {f} \in I
\rangle$. Note that this need not be a monomial ideal.
\begin{example}
Consider the polynomial $f$ in the previous example,
$f=2x+5y^2+3xyz$. Let $w=(1,1,1)$. Then initial form of
$\op{in}_w(f)=\overline{1}x+\overline{5}y^2=x+y^2$. Let $\prec$ be the
lexicographic ordering, then $\op{in}_{\prec}(\op{in}_w(f))=x$. 
\end{example}
\section{Normal form algorithm for Modules}
Let $K$ be a field with valuation and $M$ be free module over
$K[x_1,\ldots ,x_n]$ of rank $d$. Then $ M \cong \PolRing^d  $ and hence we
work only with $\PolRing^d$.

Every element of the $\PolRing^d$ can be written as
\begin{displaymath}
	\sum \limits_{k=1}^{d} \sum_{u \in \mathbb Z_{\geq0}^n} c_{u,k}x^u{e}_k\enspace,
\end{displaymath}
where $\{{e}_1,\ldots, {e}_d\}$ is the standard basis of $\PolRing^d$. 
In the above representation, $c_{u,k}x^u {e}_k$ is called a term, $x^u
{e}_k$ is called a monomial and $|u|$ is the 
degree of the term or monomial. In this case, we define  support of $f
\in \PolRing^d$, $\operatorname{supp}(f)$, to be the set $\{(u,k)
\in Z_{\geq0}^n\times [d]: c_{u,k} \neq 0\}$.  
\begin{definition}
  An element ${f} \in \PolRing^d$ is called an 
  homogeneous if every monomial occurring in ${f}$ is of same
  degree.
\end{definition}

Now, we define initial form for the elements of $\PolRing^d$.	
\begin{definition}
Let $f= \sum\limits_{(u,k) \in \Lambda} c_{u,k}x^ue_k \in \PolRing^d
$, where $\Lambda$ is the support of $f$. Let $w \in 
\Gamma^n$, $v$ be a valuation of $K$ and   $W= \min\limits_{(u,k) \in
  \Lambda} \{v(c_{u,k})+w.u\}$, where $(u,k) \in \Lambda$. The initial
form with respect to to $w$ is defined as  
\begin{displaymath}
\op{in}_w(f) = \sum_{ \substack{v(c_{u,k})+w.u=W \\ 
(u,k) \in \Lambda} }  
\overline{c_{u,k}t^{-v(c_{u,k})}}x^ue_k.
\end{displaymath}
\end{definition}

We fix an ordering on the standard bases as 
${e}_1 \prec {e}_2 \prec \ldots \prec {e}_m$ and $\prec$ be a monomial order.
 We say that $ x^{\alpha}
{e}_i \prec  x^{\beta}  {e}_j$ if $
x^{\alpha} {\prec}  x^{\beta}$ or if $x^{\alpha} =
 x^{\beta}$ and ${e}_i \prec {e}_j$. 

Let $x^ue_k$ be the monomial  in $\op{in}_\prec (\op{in}_w(f))$. Then
$\op{in}_{w,\prec}(f)$ represents the term $c_{u,k}x^ue_k$. 

\begin{example}
Let $f=2x^3e_1+12xye_2$ and $w=(0,0)$. Then with $2$-adic valuation,
$\op{in}_w(f)=x^3e_1 $,  $\op{in}_\prec(\op{in}_w(f))=x^3e_1 $ and
$\op{in}_{w,\prec}(f)=2x^3e_1$. 
\end{example}

Now having the definition for initial form of ${f}$, we 
can define the initial submodule for a submodule of $\PolRing^d$. 

\begin{definition}
Let $I$ be a submodule of $\PolRing^d$. The initial submodule of
$I$, $\op{in}_w(I)$ is defined as a submodule generated by the initial
forms of the elements of $I$, i.e $\langle \op{in}_w({f}): {f} \in I
\rangle$.   
\end{definition} 
Now we define the Gr\"{o}bner basis for a submodule
\begin{definition}
  Let $I$ be submodule of $\PolRing^d$. A generating set for $I$,
  $G= \{{g}_1, \ldots, {g}_n \}$ is a Gr\"obner basis for $I$ iff $
  \op{in}_w(I)  = \langle \op{in}_w({g}_1), \ldots, \op{in}_w({g}_n)
  \rangle$. 
  \end{definition}
\begin{lemma}
 Let $I$ be a submodule of $\PolRing^d$. Let $\prec$ be a monomial
 order. If $\{g_1, \ldots, g_s\}$ is a generating set for $I$ such
 that  $\{\op{in}_w(g_1),\ldots,\op{in}_w(g_s)\}$ is a Gr\"{o}bner
 basis of $\op{in}_w(I)$ with respect to $\prec$, then
 $\{g_1,\ldots,g_s\}$ is a Gr\"{o}bner basis of $I$ with respect to
 $w$. 
 
\end{lemma}
\proof
Since, $\{\op{in}_w(g_1),\ldots,\op{in}_w(g_s)\}$ is a  Gr\"{o}bner
basis of $\op{in}_w(I)$, it is also a generating set. Therefore by
Definition 3.5, it is a Gr\"{o}bner basis for $I$ with respect to $w$.  
\endproof

\begin{definition}
Let  $c_{\alpha,i} x^\alpha {e}_i$, $c_{\beta,j} x^\beta {e}_j$ be
terms in $\PolRing^d$, we say  $c_\beta x^\beta {e}_j$ divides
$c_\alpha x^\alpha {e}_i$ if $i=j$ and $c_\beta x^\beta$  divides
$c_\alpha x^\alpha$ 
\end{definition}

Now, we are ready to present the division algorithm. For the division
algorithm, we will need a notion of ecart function similar to the
tangent cone algorithm. 

\begin{definition}
  Let ${f},{g} \in \PolRing^d$  then 
  $\operatorname{ecart}({f},{g})= |\operatorname{supp}({g})- \operatorname{supp}({f})|.$
  \end{definition}  

\begin{theorem}
  Let ${f} \in \PolRing^d$ be a homogeneous element and $S= \{
  {g}_1, \ldots, {g}_s \}$ be a set of homogeneous elements of
  $\PolRing^d$. Then Algorithm~1  computes ${r}$ and $h_1,\ldots,
  h_s \in K[x_1, \ldots, x_n]$ such that 
\begin{displaymath}  
{f}= \sum\limits_{i=1}^{s} h_i {g}_i + {r} ,
\end{displaymath}
where $\op{in}_{w,\prec}({r}) \geq \op{in}_{w,\prec}({f})$,
$\op{in}_{w,\prec}(h_i{g}_i) \geq \op{in}_{w,\prec}({f})$ and no
monomial of  
 ${r}$ is  divisible by $\op{in}_{w,\prec}({g}_i)$ for $i \in \{1, \ldots, s \}$.
  \end{theorem}

\begin{algorithm}
\caption{Division algorithm for modules}
\label{Moduledivision}
\begin{algorithmic}[1]
\State \textbf{Input}: A finite set $B$ of homogeneous elements of $\PolRing^d$.
\State \textbf{Output} ${r}$ as mentioned Theorem 3.7.
\State Initialize:  $D = \{ {g}_1, \ldots, {g}_s \}$, $h_1 = \ldots =
h_s = 0$, ${q} = {f}$,  ${r}= 0$ 

\While {${q} \neq 0$} 

\If {there is no ${g} \in D$ with $\op{in}_{w,\prec}({g})$ dividing $\op{in}_{w,\prec}({q})$}
\State $D = D \cup \{{q}\}$
\State ${r}  = {r} +\op{in}_{w,\prec}({q})$,  ${q}= {q}-\op{in}_{w,\prec}({q})$ 
\Else
\State Choose ${g} \in D$ such that $\op{in}_{w,\prec}({g})$ divides $\op{in}_{w,\prec}({q})$ with $\operatorname{ecart}({g},{q})$ minimal.

\If {$\operatorname{ecart}({g},{q})>0$}
\State $D = D \cup \{{q}\}$
\EndIf
\State $c= \op{in}_{w,\prec}({q})/\op{in}_{w,\prec}({g})$. ${l}={q}- c{g}$
\If{${g}={g}_k$ for some $k \in \{1, \ldots ,s\}$ }
\State ${q}={l}$
\State $h_k=h_k +c$, $h_i=h_i$ for all $i \neq k$
\State ${r}  = {r}$
\Else
\State ${q}=1/(1-c){l}$
\State $h_k=1/(1-c)(h_k-ch_{k,m})$
\State ${r}= 1/(1-c)( {r}-c {r}_m)$
\EndIf
\EndIf
\EndWhile
\Return ${r}$
\end{algorithmic}
\end{algorithm}
\proof
Let ${q}_j,h_{i,j},{r}_j$ represent the value of ${q},h_{i},{r}$ at the $j$\textsuperscript{th} iteration. We use induction to prove that following conditions are true for every iteration
\begin{description}
\item[C1] ${f}={q}_j + \sum\limits_{i=1}^{s} h_{i,j} {g}_i + {r}_{j}$,
\item[C2] $ \op{in}_{w,\prec}(h_{i,j}{g}_i) \geq \op{in}_{w,\prec}(f)$,
\item[C3] No term of ${r}_j$ is divisible by $\op{in}_{w,\prec}({g}_i)$, and
\item[C4] $\op{in}_{w,\prec}({r}_j) \geq \op{in}_{w,\prec}({f})$.
\end{description}

Before the beginning of the while loop, the above conditions are
satisfied. The proof then follows by induction. 
\lspace
If there is no ${g} \in D$ with $\op{in}_{w,\prec}({g})$ dividing $\op{in}_{w,\prec}({q}_j)$, then 
${r}_j+{q}_j= {r}_{j+1} +{q}_{j+1}$. Since. the C1 holds true at $j$\textsuperscript{th} 
iteration, it is also true for $j+1$\textsuperscript{th} iteration. C2 holds true true because there is no change in 
$h_{i,j}$. C3 is true because the new monomial which is being added to ${r}_j$ is not 
divisible by $\op{in}_{w,\prec}({g}_i)$. C4 is true by the the definition of initial form.
\lspace

Suppose there exists ${g} \in D$ with $\op{in}_{w,\prec}({g})$ dividing $\op{in}_{w,\prec}({q}_j)$ and ${g}=
{g}_k$ for some $k \in \{1, \ldots, s\}$. Since, ${q}_j + h_{k,j}{g}_k= {q}_{j
+1} + h_{k,j+1}{g}_k$, C1 is satisfied. As ${r}_j$ does not change, the C2 
and 3 are satisfied. Since $\op{in}_{w,\prec}({q}_j) \geq \op{in}_{w,\prec}({f}) $ and $\op{in}_{w,\prec}(h_{k,j}{g}_k) \geq 
\op{in}_{w,\prec}({f})$, C4 is also satisfied.
\lspace
Now suppose that the last else statement is executed during the $j+1$\textsuperscript{th} iteration. Note that ${q}_j$
is a homogeneous element of $M$ and its degree remains same at all iterations. This implies $c$ is a constant term. Consider the the following equations

$${f}={q}_j +\sum\limits_{i=1}^{s} h_{i,j} {g}_i + {r}_{j} ,$$
and
$${f}={q}_m +\sum\limits_{i=1}^{s} h_{i,m} {g}_i + {r}_{m}.$$

Multiplying the second equation by $c$ and then subtracting it from the first equation we see that 
C1 is satisfied. Now let $\op{in}_{w,\prec}({q}_m)= c_mx^\alpha {e}_k$ and $\op{in}_{w,\prec}({q}_j)= 
c_jx^\alpha {e}_k$.  Since $\operatorname{val}(c_j)+w.\alpha > \operatorname{val}(c_m)+w.\alpha$. We get $\operatorname{val}(c) >  0$, since 
$c=c_j/c_m$. This implies $\operatorname{val}(1/(1-c))=0$. Since, $\operatorname{val}(c)>0$ and $\operatorname{val}(1-c)=0$, we can see that the 
C4 is satisfied. Now, since no term of ${r}_j$ and ${r}_m$ is divisible by 
$\op{in}_{w,\prec}({g}_i)$, C3 is satisfied.
\lspace
Now let $s({q}_j)$ denote the set of non-zero monomials of
$s({q}_j)$. Now, since
${q}_j$ is homogeneous polynomial there are only finitely many values
for  $s({q}_j)$. So by pigeonhole principle there exists a $j$ such
that after $j$\textsuperscript{th}  iteration the values of $s({q}_j)$
will be from a fixed set of  monomials. So, there will be $j' < j$
such that $s({q}_j)=s({q}_{j'})$ and therefore
$\operatorname{ecart}({q}_{j'},{q}_j)=0$. Now, $s({q}_{j+1})
\subsetneq s({q}_{j})$, since  the leading term is removed from
${q}_j$. So, the algorithm terminates. 
\endproof

\begin{example}
Consider $\mathbb {Q} [x,y]^2$ with $2$-adic valuation and $w= (1,1)$ and lex ordering $\succ$. Let ${f}= \left[ \begin{array}{c} 5x^3\\ 7y^3 \\ \end{array} \right] $ and ${g}_1= \left[ \begin{array}{c} 2x^2\\ 3y^2 \\ \end{array} \right] $,  ${g}_2= \left[ \begin{array}{c} 2x\\ 5y \\ \end{array} \right] $ and $D=\{{g}_1,{g}_2\}$. We can write ${f}=5x^3{e}_1+7y^3{e}_2$, ${g}_1=2x^2{e}_1+3y^2{e}_2$, ${g}_2=2x{e}_1+5y{e}_2$. Now let us calculate  $\op{in}_{w,\prec}({f})$. We have $\operatorname{val}(5)+ (1,1).(3,0) = \operatorname{val}(7)+ (1,1).(0,3)$. But $x^3 \succ y^3$. We get $\op{in}_{w,\prec}({f})= 5x^3{e}_1$.
Similarly we get $\op{in}_{w,\prec}({g}_1)=2x^2 {e}_1$ and $\op{in}_{w,\prec}({g}_2)=2x {e}_1$.
Let ${q}_0={f}=5x^3{e}_1+7y^3{e}_2$ and ${r}_0=0$. Now $\op{in}_{w,\prec}({g}_1)$ divides $\op{in}_{w,\prec}({f})$. So, we get ${q}_1={q}_0-2.5x{g}_1= 7y^3{e}_2-7.5xy^2{e}_2$ and ${r}_1=0$. $D=D \cup \{{q}_0\}$. Since there exists no ${g} \in D$ such that $\op{in}_{w,\prec}({g})$ divides $\op{in}_{w,\prec}({q}_1)$, we get ${{q}_2}= -7.5xy^2{e}_2$ and $r=7y^3{e}_2$.  $D=D \cup \{{q}_1\}$. Since there exists no ${g} \in D$ such that $\op{in}_{w,\prec}({g})$ divides $\op{in}_{w,\prec}({q}_1)$, we get ${{q}_3}= 0$ and   ${r}_3=7y^3{e}_2-7.5xy^2{e}_2$. Therefore, ${r}=7y^3{e}_2-7.5xy^2{e}_2$.

\end{example}
\section{Computation of Gr\"obner basis for submodules}
\begin{definition}
Let $c_\alpha x^\alpha {e}_i$, $c_\beta x^\beta {e}_j$ be monomials in
$\PolRing^d$. If $i=j$, then we define $\op{LCM}(c_\alpha x^\alpha
{e}_i,c_\beta x^\beta {e}_j )=\op{LCM}(x^\alpha,x^\beta){e}_j$
otherwise LCM is $0$.  
\end{definition}

Similar to S-polynomials, we define S-form for any two elements in $\PolRing^d$. 
\begin{definition}
Let ${f}, {g}$ be two elements of $\PolRing^d$. Let $x^\alpha
e_j=\op{LCM}(\op{in}_{w,\prec}({f}), \op{in}_{w,\prec}({g}))$. Then
S-form of ${f}, {g}$, is given by
\begin{displaymath}
\op{S-form}({f}, {g})=\frac{x^\alpha
  e_j}{\op{in}_{w,\prec}({f})}{f}-\frac{x^\alpha
  e_j}{\op{in}_{w,\prec}({g})}{g}.
\end{displaymath}
\end{definition}

\begin{theorem}
Let $V$ be an $n$ dimensional vector space over $K$. Let $v_1, \ldots
v_s \in V$ and $c \in K^s$. Consider the polynomial
$f_{{c}}= \sum\limits_{i=1}^{s} c_ix_i$ and let
$\op{trop}(f_{{c}})$ represents its tropicalization. Then for
every $v$ in the subspace generated by  $v_1, \ldots ,v_s $ and $w \in
\mathbb{R}^s$, there exists a  ${c} \in K^s$ with
$\sum\limits_{i=1}^{s} c_iv_i=v$ such the the value of function
$\op{trop}(f_{{c}})(w)$ is maximized. 
\end{theorem}
\proof Consider a ${c} \in K^s$ such that
$\sum\limits_{i=1}^{s} c_iv_i=v$. Assume that $v_1, \ldots ,v_s $ are
linearly dependent, otherwise the proof is trivial. Let ${c}'
\in K^s$ such that $\sum\limits_{i=1}^{s} c'_iv_i=0$. Relabel the
vectors such that $\operatorname{val}(c'_1)+w_1=
\op{trop}(f_{{c}})(w)$. Now, there exists an $\lambda$ such
that $c_1=\lambda c'_1$. From this, we get
\begin{displaymath}
\op{trop}(f_{{c}-\lambda{c'} })(w) \geq
\op{trop}(f_{{c}})(w)
\end{displaymath}
Since, ${c}-\lambda {c'}$ has less no zero components than
${c}$, we ultimately get ${b}$ such that
$\sum\limits_{i=1}^{s} b_iv_i=v$ and $v_i$ with non-zero $b_i$'s are
linearly independent. Since, the set $\{v_1, \ldots, v_s\}$ can have
only finitely many linearly independent subsets, the theorem follows. 
\endproof

Now, we state the Buchberger-like criterion for the Gr\"{o}bner basis
of submodules of $\PolRing^d$. 

\begin{theorem}
Let $S=\{g_1,\ldots,g_s\}  \subset \PolRing^d$ be finite subset and let $I$ be the
submodule generated by $S$ in $\PolRing^d$. 
If the remainder of $\op{S-form}({g}_i,{g}_j)$ on division by $S$ is
$0$ for all ${g}_i, {g}_j \in S$ then $S$ is a Gr\"obner bases for the
submodule $I$. 
\end{theorem}
\proof
Suppose ${f} \in I$. Then ${f}$ can be written as $\sum\limits_{i=1}^{s} h_i{g}_i,\; \op{for} \;  g_i \in S \; \op{and} \; h_i \in K[x_1,\ldots, x_n] \; \op{for} \; i \in \; \{1,\ldots,s\}.$  Now , 
$\op{in}_{w,\prec}({f}) \geq \min\limits_{1\leq i \leq s} ( \op{in}_{w,\prec}(h_i{g}_i))$. Now let $\op{in}_{w,\prec}({f})= c_vx^v{e}_k$ and 
$\op{in}_{w,\prec}(h_i{g}_i)= c_ux^u{e}_j$. Then we get $\operatorname{val}(c_v)+w.v \geq \min\limits_{1\leq i \leq s} (\operatorname{val}(c_i)+w.u_i)$. Now 
using the previous theorem, we choose $h_i$ such that $\min\limits_{1\leq i \leq s} ( \op{in}_{w,\prec}(h_i{g}_i))$ is maximized and let it be denoted by ${m}$. Suppose $\op{in}_{w,\prec}({f})=  {m}$,then we are done. Otherwise, $\op{in}_{w,\prec}({f})> 
{m}$. Consider the set $S=\{i:\op{in}_{w,\prec}(h_i{g}_i)= {m}\}$. Now ${g}= \sum_{i\in S} \op{in}_{w,\prec}(h_i){g}_i$. But $\op{in}_{w,\prec}({g}) < {m}$, therefore 
$${g}= \sum\limits_{i,j \in S, i \neq j} c_{i,j}\op{S-form}(x_i {g}_i,x_j{g}_j) .$$
Now, ${x}=\op{LCM}(\op{in}_{w,\prec}(x_i{g}_i),\op{in}_{w,\prec}(x_j{g}_j))$
So, $$\op{S-form}(x_i {g}_i,x_j{g}_j)= \frac{{x}}{\op{in}_{w,\prec}(x_i{g}_i)}x_i{g}_i-\frac{{x}}{\op{in}_{w,\prec}(x_j{g}_j)}x_j{g}_j$$
$$=\frac{{x}}{{x}_{ij}}\op{S-form}({g}_i,{g}_j)$$
where ${x}_{i,j}= \op{LCM}(\op{in}_{w,\prec}({g})_i,\op{in}_{w,\prec}({g})_j)$. Since $\op{S-form}({g}_i,{g}_j)$ reduces to $0$, it can be written as a sum of ${g}_i, \ldots, {g}_j$. Substituting this into ${g}$ and ${f}$, we get a representation  of ${f}$ as $\sum\limits_{i=1}^{s} h'_i{g}_i$ but $\min\limits_{1\leq i \leq s}( \op{in}_{w,\prec}(h'_i {g}_i))> {m}$, a contradiction.  
\endproof

The above criterion gives us the following algorithm for computing the Gr\"obner basis.

\begin{algorithm}[H]
\caption{Algorithm for Gr\"obner basis for modules}
\label{GrobnerBasisAlgo}
\begin{algorithmic}[1]
\State {Input}: A finite set $B$ generating the submodule $I$ of $\PolRing^d$
\State {Output} A Gr\"obner basis for the submodule $I$
\State Initialize $G=B$
\State Initialize $C=G\times G$
\While {$C \neq \phi$}
\State Choose a pair $({f},{g})$ from $C$
\State $C:=C - \{({f},{g})\}$
\State Divide $\op{S-form}({f},{g})$ by $G$ using the Algorithm 1. Let the remainder be ${r}$
\If {${h} \neq 0$} 
\State $C:=C\cup G \times \{{h}\}$
\State $G:=G \cup \{{h}\}$
\EndIf
\EndWhile
\Return $G$
\end{algorithmic}
\end{algorithm}
\proof
To prove this, we use the ascending chain condition on the module
$\PolRing^d$. Let $G_i$ represent the set $G$ in the algorithm at
$i$\textsuperscript{th} iteration. As the algorithm progresses we get
the following strictly increasing set of elements in $\PolRing^d$. 
\begin{displaymath}
G_1 \subsetneq  G_2 \subsetneq \ldots
\end{displaymath}

Let $G_i= G_{i-1} \cup \{r\}$. By Algorithm 1,
$\op{in}_{w,\prec}({r})$ is not divisible by the initial form of any
of element in $G_i$. Let $\op{in}_{w,\prec}(G)=\langle
\op{in}_{w,\prec}({g}):{g} \in G\rangle$. Therefore,
$\op{in}_{w,\prec}(G_i)  \subsetneq \op{in}_{w,\prec}(G_{i+1}) $. So,
we get an ascending chain of submodules,
\begin{displaymath}
\op{in}_{w,\prec}(G_1)  \subsetneq  \op{in}_{w,\prec}(G_2)  \ldots.
\end{displaymath}
By noetherian condition, this chain must stabilize at one point. Once
the algorithm terminates, $\op{S-form}({f},{g})$ for any ${f},{g} \in
G$ reduces to zero on division by $G$. So, by previous theorem it is a
Gr\"obner basis. 
\endproof
\begin{remark}
One motivation for studying Gr\"obner basis for fields with valuation
is that they can lead to a smaller Gr\"obner basis. Consider the
module $\mathbb{Q}[x,y,z]^r$, where $r\geq 2$. Consider the submodule
$I=\langle f,g, h \rangle$, such that $f,g$ and $h$ are of degree
$2\epsilon$. Every coefficient of $f$ except $x_1^\epsilon x_2^\epsilon e_1$ has   positive
$2$-adic valuation. Every coefficient of $g$ except $x_2^\epsilon x_3^\epsilon e_2$
has positive $2$-adic valuation and every coefficient of $h$ except
$x_1^\epsilon x_3^\epsilon e_3$ has positive $2$-adic valuation. Let $w=(0,0,0)$. Then
the initial ideal, $\op{in}_w(I)=\langle x_1^\epsilon  x_2^\epsilon e_1, x_2^\epsilon x_3^\epsilon e_2,
x_1^\epsilon x_3^\epsilon e_3 \rangle$. So, the number of generators of $\op{in}_w(I)$
remains fixed and does not  increase with $\epsilon$. Such a bound is not
possible with standard Gr\"obner basis, the number of generators will
increase with $\epsilon$. One particular example was shown in
\citep{chan2013groebner}. 
\end{remark}

\begin{remark}
Note that the initial ideal $\op{in}_w(I)$ here is dependent on the valuation of the underlying field. So, if the valuation changes  $\op{in}_w(I)$ also changes generally. But if $\op{in}_w(I)$ does not contain a monomial and $w$ lies on the unbounded of part of the tropical variety of $I$, then $\op{in}_w(I)$ doesn't contain a monomial even if the valuation changes \citep{fink2013tropical}.
\end{remark}

\section{Computation of Hilbert polynomials}
In this section we show how to compute Hilbert polynomials of  modules
using the theory described in the previous section. The standard
strategy for computing Hilbert function is to reduce to the case of
monomial ideals. We saw in Remark 4.5 that  Gr\"obner basis in the case of fields
with valuation can lead to very small monomial ideals. One can exploit
this fact to compute the Hilbert function. In Section 3, our
initial submodule was a submodule in the free module over
$\Bbbk[x_1,\ldots,x_n]$. In this section, we take it as a
submodule in the free module over $K[x_1,\ldots,x_n]$. So, the initial
module is a submodule of $\PolRing^d$. Let $I$ be the module of
$\PolRing^d$, then let $\op{in}_{w,\prec}(I)$ denotes the submodule
$\langle \op{in}_{w,\prec}({f}): {f} \in I \rangle$ in $\PolRing^d$.

\begin{theorem}
Let $I$ be a submodule generated by homogeneous elements of
$\PolRing^d$. Let $B$ be the set of monomials that do not  appear in
$\op{in}_{w,\prec}(I)$, then the residue class elements of $B$ form a $K$-vector space basis for $\PolRing^d/I$. 
\end{theorem}
\proof

We first show that elements of $B$ are linearly independent. Suppose
they are not, then there exists $b_i \in B$, $0 \neq c_i \in K$ such
that $f= \sum \limits_{i=1}^{n} c_ib_i \in I$.

Since $f \in I$, $\op{in}_{w,\prec}(f) \in \op{in}_{w,\prec}(I)$, this implies
one of the $b_i$ is in $\op{in}_{w,\prec}(I)$, which is a
contradiction.  

To show that $B$ spans the vector space,  $\{g_1,\ldots,g_s\}$ be
a Gr\"{o}bner basis for the submodule $I$. Let $f+I$ be an element of
$\PolRing^d/I$. Let $f_\delta$ denote the homogeneous component of $f$ of degree
$\delta$. Divide $f_\delta$ by $g_1,\ldots,g_s$ using  the normal form
Algorithm~1, let $r_\delta$ be the remainder. By the property of the normal
form algorithm none of the monomials appearing in $r_\delta$ is divisible by
$\op{in}_{w,\prec}(g_i)$ for $i \in \{1, \ldots, s\}$. Since
$\{g_1,\ldots,g_s\}$ is a Gr\"{o}bner basis, this implies all
monomials of $r_\delta$ belong to $B$. Since this is true for any $\delta$, $f+I$ can
be written as a sum of  residue class elements of $B$. Therefore, they generate
$\PolRing^d/I$. 
\endproof

\begin{theorem}
Let $\PolRing^d/I$ be a finitely generated graded module, where $I$
is generated by homogeneous elements. Then the Hilbert function of
$\PolRing^d/I$ and $\PolRing^d/\op{in}_{w,\prec}(I)$ are the same. 
\end{theorem}
\proof
Let $B$ be the set of monomials not in
$\op{in}_{w,\prec}(I)$. Let  $B_\delta$, $I_\delta$ and $\PolRing^d_\delta$ represent the set
of elements of degree $\delta$ in $B, I$ and $\PolRing^d$. Now,
\begin{displaymath}
\PolRing^d/I= \bigoplus_{\delta \in \mathbb{N}} \PolRing^d_\delta/I_\delta.
\end{displaymath}

By previous theorem, residue class elements of  $B$ is a basis for $\PolRing^d/I$, so $B_\delta$
forms a basis 
for  $\PolRing^d_\delta/I_\delta$. So, $\op{dim}_k \PolRing^d_\delta/I_\delta=
|B_\delta|$. Now, the since
$\op{in}_{w,\prec}(\op{in}_{w,\prec}(I))=\op{in}_{w,\prec}(I)$, the
same argument will hold  for $\PolRing^d/\op{in}_{w,\prec}(I)$.    
\endproof

So, we have reduced the problem of computing the Hilbert polynomial of
$\PolRing^d/I$ to the problem of computing the Hilbert polynomial for
$\PolRing^d/\op{in}_{w,\prec}(I)$. The Hilbert polynomial of the quotient  
of a free module by a monomial submodule can computed by standard
methods.

\section{Gr\"{o}bner basis for modules over $\mathbb Z_{p^{\ell}}[x_1,\ldots,x_n]$}
In this section, we extend the above study to free modules over $ \mathbb Z_{p^{\ell}}[x_1,\ldots,x_n]$, where $p$ is a prime and $\ell$ is a positive integer. Let $M$ be a free module of rank $d$ over  $\mathbb Z_{p^{\ell}}[x_1,\ldots,x_n]$, then $M \cong  \mathbb Z_{p^{\ell}}[x_1,\ldots,x_n]^d$. In order  to define an ordering on the monomials of $\mathbb Z_{p^{\ell}}[x_1,\ldots,x_n]^d$, we first define the following map

\begin{definition} Let $v: \mathbb Z_{p^{\ell}} \rightarrow  \{0,1,\ldots,\ell-1, \infty\}$ be such that if $m=p^ka$ and $\op{gcd}(p,a)=1$ then $v(m)=k$ and $v(0)= \infty$.
\end{definition}

\begin{theorem} 
Let $a, b \in  \mathbb Z_{p^{\ell}}$ and $v$ is defined as Definition 6.1
\begin{enumerate}
\item $v(ab)=v(a)+v(b)$ when $ab \neq 0$, and
\item $v(a+b) \geq \operatorname{min}\{v(a), v(b)\}$ when $a+b \neq 0$.
\end{enumerate}

\end{theorem}
\proof
Let $a=p^kc, b=p^{k'}c'$, $ \op{gcd}(p,c)= \op{gcd}(p,c')=\op{gcd}(p,cc')=1$ and $k'+k <\ell$. Now, $ab=p^{k+k'}cc'=qp^\ell+r$. We get $r=p^{k+k'}(cc'-qp^{\ell-k-k'})$.  Now, $\op{gcd}(cc'-qp^{\ell-k-k'},p)=1$. So, we get $v(ab)=k+l=v(a)+v(b)$. Now, let  $k' <k$.
Then $a+b=p^kc+p^{k'}c'=qp^\ell+r$. We get $r=p^{k'}c'+p^kc-qp^\ell=p^{k'}(c'+p^{k-k'}c-qp^{\ell-l})$.  Now, $\op{gcd}(p,c'+p^{k-k'}c-qp^{\ell-l})=\op{gcd}(p,c')=1$. So, $v(a+b)=k'$.
\endproof
Now using the  map $v$ from Definition 6.1 we can define ordering on the terms in $\mathbb Z_{p^{\ell}}[x_1,\ldots,x_n]^d$.

\begin{definition}
Let $\succ$ be a monomial order for  $\mathbb Z_{p^{\ell}}[x_1,\ldots,x_n]$ and let $c_{u,i}x^u {e}_i$ and $c_{v,j}x^v {e}_j$ be two monomials. We say that $c_{u,i}x^u {e}_i \succ c_{v,j}x^v {e}_j $ if $x^u \succ x^v$ or if $x^u=x^v$ and ${e}_i \succ {e}_j $
\end{definition}

\begin{definition}
Let $w$ be a weight vector. We say that $c_{u,i}x^u {e}_i < c_{v,j}x^v {e}_j $ if $v(c_{u,i})+w.u< v(c_{v,j})+ w.v $ or  $v(c_{u,i})+w.u=v(c_{v,j})+ w.v $ and $c_{u,i}x^u {e}_i \succ c_{v,j}x^v {e}_j $
\end{definition}

\begin{definition}

Let	${f} \in \mathbb Z_{p^{\ell}}[x_1,\ldots,x_n]^d$ and  $\sum \limits_{k=1}^{d} \sum_{u \in \mathbb Z_{\geq0}^n} c_{u,k}x^u{e}_k\enspace$, $\op{in}_{w,\prec}({f})$ is $c_{u,i}x^u {e}_i$, such that $c_{u,j}x^u {e}_j \leq c_{v,i}x^v {e}_i$ for all $(v,i) \in \op{supp}(f)$.
\end{definition}

Now we define divisibility  on terms in $\mathbb Z_{p^{\ell}}[x_1,\ldots,x_n]^d$.
\begin{definition}
Let $c_\alpha x^\alpha {e}_i$, $c_\beta x^\beta {e}_j$ be monomials in $\mathbb Z_{p^{\ell}}[x_1,\ldots,x_n]^d$, we say $c_\beta x^\beta {e}_j$ divides $c_\alpha x^\alpha {e}_i$ if $i=j$ , $c_\beta$ divides $ c_\alpha  $  and  $x^\beta$ divides $x^\alpha$.
\end{definition}

In this section, we present a normal formal algorithm according to the order given mentioned in Definition 6.4. 

\begin{theorem}
  Let ${f} \in \mathbb Z_{p^{\ell}}[x_1,\ldots,x_n]^d $ be a homogeneous elements and $S= \{ {g}_1, \ldots, {g}_s \}$ be set of homogeneous elements of $\mathbb Z_{p^{\ell}}[x_1,\ldots,x_n]^d$. Then Algorithm 1  computes ${r}$ and $h_1,\ldots, h_s \in \mathbb Z_{p^{\ell}}[x_1, \ldots, x_n]$ such that
  
  $${f}= \sum\limits_{i=1}^{s} h_i {g}_i + {r} ,$$
where $\op{in}_{w,\prec}({r}) \geq \op{in}_{w,\prec}({f})$, $\op{in}_{w,\prec}(h_i{g}_i) \geq \op{in}_{w,\prec}({f})$ and no monomial of ${r}$ is  divisible by $\op{in}_{w,\prec}({g}_i)$ for $i \in \{1, \ldots, s \}$.
 \end{theorem}

\proof
Let ${q}_j,h_{i,j},{r}_j$ represent the value of ${q},h_{i},{r}$ at the $j$ \textsuperscript{th} iteration.
\lspace
The proof is similar to the proof of Theorem 3.9. The difference lies in the third if statement.
Assume that $j$\textsuperscript{th} statement is being executed. Note that ${q}_j$ is a homogeneous element of $\mathbb Z_{p^{\ell}}[x_1,\ldots,x_n]^d$ and its degree remains same at iteration. Now let  ${g}_j \in D$ be equal to ${q}_m$ for some $m$ less than $j$. Note that ${q}_j$ is a homogeneous element of $\mathbb Z_{p^{\ell}}[x_1,\ldots,x_n]^d$ and its degree remains the same at all iterations. Let $\op{in}_{w,\prec}({q}_j)=c_ux^u{e}_i$ and $\op{in}_{w,\prec}({q}_m)=c_vx^v{e}_i$.
Since, degree of ${q}_j$ and ${q}_m$ are same we get that $x^u=x^v$. Also, since $\op{in}_{w,\prec}({q}_j) < \op{in}_{w,\prec}({q}_m)$, we get $\operatorname{val}(c_u)+w.u<\operatorname{val}(c_v)+w.v$. So $\operatorname{val}(-c) > 0$. Therefore, $\op{gcd}(1-c,p)=\op{gcd}(1-c, p^\ell)=1$. This implies $1-c$ is invertible. The rest of the proof is similar to Theorem 3.9.
\endproof

\begin{example}
Let $M= \mathbb {Z}/8 \mathbb {Z} [x,y] ^2$ with  $w= ( 1,1)$ and lex ordering $\succ$. Let ${f}= \left[ \begin{array}{c} 4x^3\\ 6y^3 \\ \end{array} \right] $ and ${g}_1= \left[ \begin{array}{c} 4xy\\ 2y^2 \\ \end{array} \right] $,  ${g}_2= \left[ \begin{array}{c} 2x\\ 2y \\ \end{array} \right] $ and $D=\{{g}_1,{g}_2\}$ . We can write ${f}=4x^3{e}_1+6y^3{e}_2$, ${g}_1=4xy{e}_1+2y^2{e}_2$, ${g}_2=2x{e}_1+2y{e}_2$.
Now let us calculate  $\op{in}_{w,\prec}({f})$. We have $\operatorname{val}(6)+ (1,1).(0,3) < \operatorname{val}(4)+ (1,1).(3,0)$.  We get $\op{in}_{w,\prec}({f})= 6y^3{e}_2$.
Similarly, we get $\op{in}_{w,\prec}({g}_1)=4xy {e}_1$ and $\op{in}_{w,\prec}({g}_2)=2x {e}_1$
Let ${q}_0=4x^3{e}_1+6y^3{e}_2$ and ${r}_0=0.$ Since there exists no ${g} \in D$ such that $\op{in}_{w,\prec}({g})$ divides $\op{in}_{w,\prec}({q}_1)$, we get ${q}_1=4x^3{e}_1$ , ${r}_1=6y^3{e}_2$ and $D= D \cup \{{q}_0\}$. Now $\op{in}_{w,\prec}({g}_2)$ divides $\op{in}_{w,\prec}({q}_1)$. So, we get ${q}_2={q}_1-2x^2{g}_2= -4x^2y{e}_2$ ,${r}_2=6y^3$ and $D= D \cup \{{q}_1\}.$ Since there exists no ${g} \in D$ such that $\op{in}_{w,\prec}({g})$ divides $\op{in}_{w,\prec}({q}_2)$, we get ${q}_3=0$ and ${r}_1=6y^3{e}_2+4x^2y{e}_2$. Therefore, we get ${r}=6y^3{e}_2+4x^2y{e}_2$. 
\end{example}

\begin{definition}
Let $c_\alpha x^\alpha {e}_i$, $ c_\beta x^\beta {e}_j$ be monomials in $\mathbb Z_{p^{\ell}}[x_1,\ldots,x_n]^d$. If $i=j$ , then $\op{LCM}(c_\alpha x^\alpha {e}_i, c_\beta x^\beta {e}_j )=\op{LCM}(c_\alpha, c_\beta)(\op{LCM}(x^\alpha, x^\beta)){e}_j$ otherwise is $0$. 
\end{definition}

We define S-form of two elements of $\mathbb Z_{p^{\ell}}[x_1,\ldots,x_n]^d$ which is similarly to the definition 4.2.

Following is the  Buchberger-like criterion for Gr\"obner basis in this case

\begin{theorem}
Let $S \subset \mathbb Z_{p^{\ell}}[x_1,\ldots,x_n]^d$ be a finite subset and let $I$ be the submodule generated by $S$ in $\mathbb Z_{p^{\ell}}[x_1,\ldots,x_n]^d$. If the remainder of $\op{S-form}({g}_i,{g}_j)$ on division by $S$ is $0$ for all ${g}_i, {g}_j \in S$ then $S$ is a Gr\"obner bases for the submodule $I$.
\end{theorem}
\proof
The proof is similar to Theorem 6.4. The maximum of $\operatorname{min}( \op{in}_{w,\prec}(h_i{g}_i))$ is guaranteed because 
there are only finitely many ways to write ${f}$ as a sum of ${g}_i$ as $\mathbb Z_{p^{\ell}}$ 
is finite.
\endproof

The algorithm for computing the Gr\"obner basis  is the same as Algorithm 2. Here we present the proof of correctness.

\proof
Since $\mathbb Z_{p^{\ell}}$ is a finite ring therefore it is noetherian. Now using Hilbert basis basis theorem we get $\mathbb Z_{p^{\ell}}[x_1,\ldots,x_n]$ is noetherian. Since, $\mathbb Z_{p^{\ell}}[x_1,\ldots,x_n]^d$ is finitely generated, $\mathbb Z_{p^{\ell}}[x_1,\ldots,x_n]^d$ is a noetherian module. Rest of the goes along the same line as the proof of Algorithm 2.
\endproof
\section{Conclusion}
In this paper, we studied a generalization of Gr\"obner basis for modules that also takes the valuation of coefficients into account. We expect the algorithms presented in this paper to have many computational advantages. For example, they can lead to smaller  Gr\"obner basis.  Also, to deal with blowing up of coefficients one can first compute the Gr\"obner basis over $\mathbb Z_{p^{\ell}}$ and then lift it to the field $\mathbb{Q}$.

\bibliographystyle{jtbnew} 
\bibliography{a}

\end{document}